\documentclass[a4paper, 12pt]{article}
\usepackage{}

\usepackage{mathrsfs}
\usepackage{epsfig}
\usepackage{amsmath}
\usepackage{amssymb}
\usepackage{latexsym}
\usepackage{amsfonts}
\usepackage{amsthm}

\usepackage[numbers,sort&compress]{natbib}
\usepackage{graphicx}
\ifx\pdfoutput\undefined  \else
 \fi

\usepackage[centerlast]{caption2}

\usepackage{color}

\usepackage{txfonts}
\usepackage{amssymb}
\usepackage{mathrsfs}
\usepackage{amsfonts}

\newtheorem{theorem}{Theorem}[section]
\newtheorem{lemma}[theorem]{Lemma}

\usepackage{latexsym,amssymb}
\usepackage{amsmath}

\newcommand{\F}{{\mathbb F}}

\begin{document}

\title{Davenport constant for semigroups II}

\author{
Guoqing Wang\\
\small{Department of Mathematics, Tianjin Polytechnic University, Tianjin, 300387, P. R. China}\\
\small{Email: gqwang1979@aliyun.com}
}

\date{}
\maketitle

\begin{abstract}  Let $\mathcal{S}$ be a finite commutative semigroup. The Davenport constant of $\mathcal{S}$, denoted ${\rm D}(\mathcal{S})$, is defined to be the least positive integer $\ell$ such that every sequence $T$ of elements in $\mathcal{S}$ of length at least $\ell$ contains a proper subsequence $T'$ ($T'\neq T$) with the sum of all terms from $T'$ equaling the sum of all terms from $T$. Let $q>2$ be a prime power, and let $\F_q[x]$ be the ring of polynomials over the finite field $\F_q$.
Let $R$ be a quotient ring of $\F_q[x]$ with $0\neq R\neq \F_q[x]$. We prove that $${\rm D}(\mathcal{S}_R)={\rm D}(U(\mathcal{S}_R)),$$
where $\mathcal{S}_R$ denotes the multiplicative semigroup of the ring $R$, and $U(\mathcal{S}_R)$ denotes the group of units in $\mathcal{S}_R$.
\end{abstract}


\section {Introduction}

Let $G$ be an additive finite abelian group. A sequence $T$ of
elements in $G$ is called a {\sl zero-sum sequence} if the sum of
all terms of $T$ equals to zero, the identity element of $G$. The Davenport constant ${\rm D}(G)$ of
$G$ is defined to be the smallest positive integer $\ell$  such that,
every sequence $T$ of elements in $G$ of length at least $\ell$ contains a nonempty
subsequence $T'$ with the sum of all terms of $T'$ equaling zero. Though
attributed to H. Davenport who proposed \cite{Davenport} the study
of this constant in 1965, K. Rogers \cite{rog1} had first studied it
in 1963 and this reference was somehow missed out by most of the
authors in this area. The Davenport constant together with the celebrated Erd\H{o}s-Ginzburg-Ziv Theorem  obtained by P. Erd\H{o}s, A. Ginzburg and A. Ziv in 1961 were two pioneering
researches for Zero-sum Theory (see \cite{GaoGeroldingersurvey} for a survey) which has been developed into a branch of Combinatorial Number Theory.

 \noindent\textbf{Throrem A.}  \cite{EGZ} \ (Erd\H{o}s-Ginzburg-Ziv Theorem)
\ {\sl Every sequence of $2n-1$ elements in an additive finite
abelian group of order $n$ contains a zero-sum subsequence of length
$n$. }

During the past five decades, the Davenport constant and the Erd\H{o}s-Ginzburg-Ziv Theorem together with a large of related problems have been studied extensively for the setting of groups (see \cite{AdhikariRath,Gao96,Gao2000,GG09,GG13,GS92,Grynkiewicz} for example).
In 2008, the author of this paper and W.D. Gao formulated the definition of the Davenport constant for finite commutative semigroups which is stated as follows.

\noindent \textbf{Definition B.} \cite{wanggao} \ {\sl Let $\mathcal{S}$ be a commutative semigroup (not necessary finite). Let $T$ be a sequence of elements in $\mathcal{S}$. We call $T$  reducible if $T$ contains a proper subsequence $T'$ ($T'\neq T$) such that the sum of all terms of $T'$ equals the sum of all terms of $T$. Define the Davenport constant of the semigroup $\mathcal{S}$, denoted ${\rm D}(\mathcal{S})$, to be the smallest $\ell\in \mathbb{N}\cup\{\infty\}$ such that every sequence $T$ of length at least $\ell$ of elements in $\mathcal{S}$ is reducible.}

In fact, starting from the research of Factorization Theory in Algebra, A. Geroldinger and F. Halter-Koch in 2006 have formulated another closely related definition, ${\rm d}(\mathcal{S})$, for any commutative semigroup  $\mathcal{S}$, which is called the small Davenport constant. For the completeness, their definition is also stated here.

\noindent \textbf{Definition C.} (Definition 2.8.12 in \cite{GH}) \ {\sl For a commutative semigroup $\mathcal{S}$, let ${\rm d}(\mathcal{S})$ denote the smallest $\ell \in
\mathbb{N}_0\cup \{\infty\}$ with the following property:

For any $m\in \mathbb{N}$ and $a_1, \ldots,a_m\in \mathcal{S}$ there
exists a subset $I\subset [1,m]$ such that $|I|\leq \ell$ and
$$
\sum_{i=1}^m a_i=\sum_{i \in I}a_i.
$$}

We have the following connection between the (large) Davenport constant ${\rm D}(\mathcal{S})$ and the small Davenport constant ${\rm d}(\mathcal{S})$ for any finite commutative semigroup $\mathcal{S}$.

\noindent \textbf{Proposition D.} \  {\sl Let $\mathcal{S}$ be a finite  commutative semigroup. Then,
\begin{enumerate}
\item  $\mathsf d(\mathcal{S})<
\infty.$ (See Proposition 2.8.13 in \cite{GH}.)
\item  $\mathsf D(\mathcal{S})=\mathsf d(\mathcal{S})+1.$ (See Proposition 1.2 in \cite{AdhikariGaoWang14}.)
\end{enumerate}
}

In 2014, the author together with S.D. Adhikari and W.D. Gao \cite{AdhikariGaoWang14} also generalized the Erd\H{o}s-Ginzburg-Ziv Theorem to finite commutative semigroups.

Very recently, H.L. Wang, L.Z. Zhang, Q.H. Wang and Y.K. Qu \cite{wang-zhang-wang-qu} made a study of the Davenport constant of the multiplicative semigroup of a quotient ring of $\mathbb{F}_p[x]$. Precisely, they proved the following.

\noindent \textbf{Theorem E.} \  {\sl For any prime $p>2$, let $f(x)$ be a nonconstant polynomial of $\F_p[x]$ such that $f(x)$ factors into a product of pairwise non-associate irreducible polynomials. Let $R=\frac{F_p[x]}{(f(x))}$. Then  $${\rm D}(\mathcal{S}_R)={\rm D}(U(\mathcal{S}_R)),$$
where $\mathcal{S}_R$ denotes the multiplicative semigroup of the quotient ring $\frac{\F_p[x]}{(f(x))}$ and $U(\mathcal{S}_R)$ denotes the group of units in $\mathcal{S}_R$. }

Moreover, they conjectured that ${\rm D}(\mathcal{S}_R)={\rm D}(U(\mathcal{S}_R))$ holds true for all prime $p>2$ and any nonconstant polynomial $f(x)\in \F_p[x]$.

\medskip

In this paper, we obtained the following result for the quotient ring of the ring of polynomials over any finite field $\F_q$ where $q>2$. As a special case, we affirmed their conjecture.

\begin{theorem}\label{Theorem irreducible quotient ring}
\  Let $q>2$ be a prime power, and let $\F_q[x]$ be the ring of polynomials over the finite field $\F_q$.
Let $R$ be a quotient ring of $\F_q[x]$ with $0\neq R\neq \F_q[x]$. Then $${\rm D}(\mathcal{S}_R)={\rm D}(U(\mathcal{S}_R)),$$
where $\mathcal{S}_R$ denotes the multiplicative semigroup of the ring $R$, and $U(\mathcal{S}_R)$ denotes the group of units in $\mathcal{S}_R$.
\end{theorem}

\section{The proof of Theorem \ref{Theorem irreducible quotient ring}}

We begin this section by giving some preliminaries.

Let $\mathcal{S}$ be a finite commutative semigroup.
The operation on $\mathcal{S}$ is denoted by $+$.
The identity element of $\mathcal{S}$, denoted $0_{\mathcal{S}}$ (if exists), is the unique element $e$ of
$\mathcal{S}$ such that $e+a=a$ for every $a\in \mathcal{S}$. If $\mathcal{S}$ has an identity element $0_{\mathcal{S}}$, let
$$U(\mathcal{S})=\{a\in \mathcal{S}: a+a'=0_{\mathcal{S}} \mbox{ for some }a'\in \mathcal{S}\}$$ be the group of units
of $\mathcal{S}$. For any element $c\in\mathcal{S}$ and any subset $A\subseteq\mathcal{S}$, let $${\rm St}_{A}(c)=\{a\in A: a+c=c\}$$ denote the stabilizer of $c$ in $A$.

On a commutative semigroup $\mathcal{S}$ the Green's preorder, denoted $\leqq_{\mathcal{H}}$, is defined by
$$a \leqq_{\mathcal{H}} b\Leftrightarrow a=b \ \ \mbox{or}\ \ a=b+c$$ for some $c\in \mathcal{S}$. Green's congruence, denoted
$\mathcal{H}$, is a basic relation introduced by Green for semigroups which is defined by:
$$a \ \mathcal{H} \ b \Leftrightarrow a \ \leqq_{\mathcal{H}} \ b \mbox{ and } b \ \leqq_{\mathcal{H}} \ a.$$
For any element $a$ of $\mathcal{S}$,  let $H_a$ be the congruence class by $\mathcal{H}$ containing $a$.
We write $a<_{\mathcal{H}} b$ to mean that $a \leqq_{\mathcal{H}} b$ but $H_a\neq H_b$. The following easy fact will be used later.

\begin{lemma}\label{Lemma folklore} (folklore) \ For any element $a\in \mathcal{S}$, $U(S)$ acts on the congruence class  $H_a$ and
${\rm St}_{U(\mathcal{S})}(a)$ is a subgroup of $U(\mathcal{S})$.
\end{lemma}

In what follows, we also need some notations introduced by A. Geroldinger and F. Halter-Koch (see \cite{GH}), which are very helpful to dealing with the problems in zero-sum theory and factorization theory.

The sequence $T$ of elements in the semigroups $\mathcal{S}$ is denoted by $$T=a_1a_2\cdot\ldots\cdot a_{\ell}=\prod\limits_{a\in \mathcal{S}} a^{\ {\rm v}_a(T)},$$ where ${\rm v}_a(T)$ denotes the multiplicity of the element $a$ in the sequence $T$. By $\cdot$ we denote the operation to join sequences.
Let $T_1,T_2$ be two sequences of elements in the semigroups $\mathcal{S}$. We call $T_2$
a subsequence of $T_1$ if $${\rm v}_a(T_2)\leq {\rm v}_a(T_1)$$ for every element $a\in \mathcal{S}$, denoted by $$T_2\mid T_1.$$ In particular, if $T_2\neq T_1$, we call $T_2$ a {\sl proper} subsequence of $T_1$, and write $$T_3=T_1  T_2^{-1}$$ to mean the unique subsequence of $T_1$ with $T_2\cdot T_3=T_1$.  Let $$\sigma(T)=a_1+a_2+\cdots+a_{\ell}$$ be the sum of all terms in the sequence $T$.
By $\lambda$ we denote the
empty sequence.
If $S$ has an identity element $0_{\mathcal{S}}$,  we allow $T=\lambda$ and adopt the convention
that $\sigma(\lambda)=0_\mathcal{S}$.
We say that $T$ is {\it
reducible} if $\sigma(T')=\sigma(T)$ for some proper subsequence $T'$ of $T$
(note that, $T'$ is probably the empty sequence $\lambda$ if $\mathcal{S}$
has the identity element $0_{\mathcal{S}}$ and $\sigma(T)=0_{\mathcal{S}}$). Otherwise, we call $T$
{\it irreducible}. For more related terminology used in additive problems for semigroups, one is refereed to  \cite{AdhikariGaoWang14,wang}. Here, the following two lemmas are necessary.

\begin{lemma}(\cite{GH}, Lemma 6.1.3) \label{Lemma recusive Davenport constant} \ Let $G$ be a finite abelian group, and let $H$ be a subgroup of $G$. Then, ${\rm D}(G)\geq {\rm D}(G/H)+{\rm D}(H)-1$.
\end{lemma}

\begin{lemma}\label{proposition D(U(G))leq D(G)} (see \cite{wanggao}, Proposition 1.2) \
Let $\mathcal{S}$ be a finite commutative semigroup with an identity. Then ${\rm D}(U(\mathcal{S}))\leq
{\rm D}(\mathcal{S})$.
\end{lemma}

\bigskip

Now we are in a position to  prove Theorem \ref{Theorem irreducible quotient ring}.

\noindent {\sl Proof of Theorem  \ref{Theorem irreducible quotient ring}.} \ By Lemma \ref{proposition D(U(G))leq D(G)}, we need only to show that
$${\rm D}(\mathcal{S}_R)\leq {\rm D}(U(\mathcal{S}_R)).$$

Since the ring $\mathbb{F}_q[x]$ is a principal ideal domain and $0\neq R\neq \mathbb{F}_q[x]$, we have that $R=\mathbb{F}_q[x]\diagup (f)$ for some nonconstant {\sl monic} polynomial $f\in \mathbb{F}_q[x]$. Let
\begin{equation}\label{equation factorization of f(x)}
f=f_1^{n_1}*f_2^{n_2}*\cdots * f_r^{n_r}
\end{equation}
be the factorization of $f(x)$ in $\mathbb{F}_q[x]$, where $r\geq 1$, $n_1,n_2,\ldots,n_r\geq 1$, and
$f_1, f_2, \ldots,f_r$ are pairwise non-associate monic irreducible polynomials of $\mathbb{F}_q[x]$. To proceed, we need to introduce some notations.

Take an arbitrary element $a\in \mathcal{S}_{R}$.
Let $\theta_a\in \mathbb{F}_q[x]$ be the unique polynomial corresponding to the element $a$ with the least degree, i.e.,
$$\overline{\theta_a}=\theta_a+(f)$$ is the corresponding form of $a$ in the quotient ring $R$ with $${\rm deg}(\theta_a)\leq {\rm deg}(f)-1.$$
By $\gcd(\theta_a,f)$ we denote the greatest common divisor  of the two polynomials $\theta_a$ and $f$ in $\mathbb{F}_q[x]$ ({\sl the unique monic polynomial with the greatest degree which divides both $\theta_a$ and $f$}), in particular, by \eqref{equation factorization of f(x)},  $$\gcd(\theta_a,f)=f_1^{\alpha_1}*f_2^{\alpha_2}*\cdots * f_r^{\alpha_r}$$ where $\alpha_i\in [0,n_i]$ for $i=1,2,\ldots,r$.

\noindent $\bullet$ {\small For notational convenience,  we shall write ${\rm St}_{U(\mathcal{S}_R)}(\cdot)$ simply as ${\rm St}(\cdot)$ in what follows.
It is also noteworthy that for any $a,b,c\in \mathcal{S}_R$, $a+b=c$ holds if and only if $\theta_a*\theta_b\equiv \theta_c \pmod {f}$.}

Now we prove the following claim.

\noindent \textbf{Claim A.} \ {\sl Let $a$ and $b$ be two elements of $\mathcal{S}_R$. Then the following conclusions hold:

{\rm (i)} \ If $a \leqq_{\mathcal{H}} b$ then $\gcd(\theta_b,f)\mid \gcd(\theta_a,f)$
and ${\rm St}(b)\subseteq{\rm St}(a);$

{\rm (ii)} \ $a \ {\mathcal{H}} \  b \Leftrightarrow \gcd(\theta_b,f)=\gcd(\theta_a,f) \Leftrightarrow {\rm St}(b)={\rm St}(a).$

\noindent {\sl Proof of Claim A.} \  Assume $a \leqq_{\mathcal{H}} b$.  Since $\mathcal{S}_R$ has the identity element $0_{\mathcal{S}_R}$, we have $$a=b+c \ \ \mbox{for some }c\in \mathcal{S}_R.$$ It follows that
$$\gcd(\theta_b,f)\mid \gcd(\theta_b*\theta_c,f)=\gcd(\theta_a,f).$$
For any element $d\in {\rm St}(b)$, $d+a=d+(b+c)=(d+b)+c=b+c=a$,
and so $d\in {\rm St}(a)$. It follows that $${\rm St}(b)\subseteq {\rm St}(a).$$ This proves Conclusion {\rm (i)}.

Now we prove Conclusion {\rm (ii)}.

Assume $a \ {\mathcal{H}} \  b$. Then $a \ \leqq_{\mathcal{H}} \  b$ and $b \ \leqq_{\mathcal{H}} \  a$. It follows from Conclusion {\rm (i)} that $$\gcd(\theta_b,f)=\gcd(\theta_a,f)$$ and $${\rm St}(b)={\rm St}(a).$$

Assume $\gcd(\theta_b,f)=\gcd(\theta_a,f).$
It follows that there exist polynomials $h, h'\in \mathbb{F}_q[x]$ such that $$\theta_a* h\equiv \theta_b\pmod {f}$$
and
$$\theta_b* h'\equiv \theta_a\pmod {f}.$$
It follows that $b \leqq_{\mathcal{H}}  a$ and $a \leqq_{\mathcal{H}}  b$, i.e., $$a \ \mathcal{H} \ b.$$

Assume ${\rm St}(b)={\rm St}(a)$.  To prove $a \ \mathcal{H} \ b$, we suppose to the contrary that $a \ \mathcal{H} \ b$ does not hold. Then
$\gcd(\theta_b,f)\neq \gcd(\theta_a,f).$ We may suppose without loss of generality that
there exist integers $k\in [1,r]$ and $m_k\in [1,n_k]$ such that
\begin{equation}\label{equation construction 1}
f_k^{m_k}\mid \gcd(\theta_a,f)
\end{equation} and
\begin{equation}\label{equation construction 2}
f_k^{m_k}\not\mid \gcd(\theta_b,f).
\end{equation}
Let
\begin{equation}\label{equation construction 3}
h=\frac{f}{f_k^{m_k}}.
\end{equation}
Take an element $\xi\in \mathbb{F}_q\setminus \{0_{\mathbb{F}_q}, 1_{\mathbb{F}_q}\}$.

Now we show that
\begin{equation}\label{equation h(x)+1,f(x)=1}
\gcd(h+1_{\mathbb{F}_q},f)=1_{\mathbb{F}_q}
\end{equation}
or
\begin{equation}\label{equation 2h(x)+1,f(x)=1}
 \gcd(\xi*h+1_{\mathbb{F}_q},f)=1_{\mathbb{F}_q}.
 \end{equation}
Suppose to the contrary that $\gcd(h+1_{\mathbb{F}_q},f)\neq 1_{\mathbb{F}_q}$ and $\gcd(\xi*h+1_{\mathbb{F}_q},f)\neq 1_{\mathbb{F}_q}$.
By \eqref{equation factorization of f(x)} and \eqref{equation construction 3}, we have that $f_i\not\mid \gcd(h+1_{\mathbb{F}_q},f)$ and $f_i\not\mid \gcd(\xi*h+1_{\mathbb{F}_q},f)$ for each $i\in[1,r]\setminus \{k\}$, and thus $f_k\mid (h+1_{\mathbb{F}_q})$ and $f_k\mid (\xi*h+1_{\mathbb{F}_q})$. It follows that $f_k\mid \xi*(h+1_{\mathbb{F}_q})-(\xi*h+1_{\mathbb{F}_q})=\xi-1_{\mathbb{F}_q}$, which is absurd. This proves that \eqref{equation h(x)+1,f(x)=1} or \eqref{equation 2h(x)+1,f(x)=1} holds.

Take an element $d\in \mathcal{S}_R$ with $$\theta_d\equiv h+1_{\mathbb{F}_q}\pmod {f}$$ or $$\theta_d\equiv \xi*h+1_{\mathbb{F}_q} \pmod {f}$$ according to \eqref{equation h(x)+1,f(x)=1} or \eqref{equation 2h(x)+1,f(x)=1} holds respectively. It follows that $$d\in U(\mathcal{S}_R),$$ and follows from \eqref{equation construction 1}, \eqref{equation construction 2} and \eqref{equation construction 3} that
$$\theta_a*\theta_d\equiv \theta_a\pmod {f}$$
and $$\theta_b*\theta_d\not\equiv \theta_b\pmod {f}.$$ That is, $d\in {\rm St}(a)\setminus {\rm St}(b),$ a contradiction with ${\rm St}(a)={\rm St}(b)$. Hence, we have that $$a \ \mathcal{H} \ b.$$
This proves Claim A. \qed

Let $T=a_1a_2\cdot\ldots\cdot a_{\ell}$ be an arbitrary sequence of elements in $\mathcal{S}_R$ of length $$\ell={\rm D}(U(\mathcal{S}_R)).$$ It suffices to show that $T$ contains a proper subsequence $T'$ with $\sigma(T')=\sigma(T)$.

Take a shortest subsequence $V$ of $T$ such that
\begin{equation}\label{equation sigma(V)Hsigma(T)}
\sigma(V) \ \mathcal{H} \ \sigma(T).
\end{equation}
We may assume without loss of generality that $$V=a_1\cdot a_2\cdot\ldots\cdot a_t\ \ \ \mbox{where} \ \ t\in [0,\ell].$$ By the minimality of $|V|$, we derive that $$0_{\mathcal{S}_R}>_{\mathcal{H}}a_1>_{\mathcal{H}}(a_1+a_2)>_{\mathcal{H}}>\cdots>_{\mathcal{H}}\sum_{i=1}^t a_i.$$
Denote $$K_0=\{0_{\mathcal{S}_R}\}$$
and
$$K_i={\rm St}(\sum\limits_{j=1}^i a_j) \ \ \ \mbox{for each}\ \ \ i\in [1,t].$$
By Lemma \ref{Lemma folklore}, $K_i$ is a subgroup of $U(\mathcal{S}_R)$ for each $i\in [1,t]$. Moreover, since ${\rm St}(0_{\mathcal{S}_R})=K_0$,
it follows from Claim A that
$$K_0\lneq K_1\lneq K_2\lneq \cdots\lneq K_t.$$
For $i\in [1,t]$, since $\frac{U(\mathcal{S}_R)}{K_i}\cong \frac{U(\mathcal{S}_R)\diagup K_{i-1}}{K_{i}\diagup K_{i-1}}$ and ${\rm D}(K_i\diagup K_{i-1})\geq 2$, it follows from Lemma \ref{Lemma recusive Davenport constant} that
$$\begin{array}{llll}
{\rm D}(U(\mathcal{S}_R)\diagup K_i)&=& {\rm D}(\frac{U(\mathcal{S}_R)\diagup K_{i-1}}{K_{i}\diagup K_{i-1}}) \\
&\leq & {\rm D}(U(\mathcal{S}_R)\diagup K_{i-1})-({\rm D}(K_i\diagup K_{i-1})-1)\\
&\leq & {\rm D}(U(\mathcal{S}_R)\diagup K_{i-1})-1.\\
\end{array}$$
It follows that \begin{align}\label{equation length and D()}
\begin{array}{llll}
1\leq {\rm D}(U(\mathcal{S}_R)\diagup K_t)&\leq & {\rm D}(U(\mathcal{S}_R)\diagup K_{t-1})-1 \\
& \vdots & \\
&\leq & {\rm D}(U(\mathcal{S}_R)\diagup K_0)-t\\
&=& {\rm D}(U(\mathcal{S}_R))-t\\
&=& \ell-t\\
&=& |TV^{-1}|.\\
\end{array}
\end{align}
By \eqref{equation sigma(V)Hsigma(T)} and Conclusion {\rm (ii)} of Claim A, we have
\begin{equation}\label{equation two common divisors equal}
\gcd(\theta_{\sigma(V)},f)=\gcd(\theta_{\sigma(T)},f).
\end{equation}
Let $$\mathcal{J}=\{j\in [1,r]: f_j^{n_j}\mid \theta_{\sigma(T)}\}.$$ By \eqref{equation two common divisors equal}, we have that
\begin{equation}\label{equation fi(x)notmid a}
f_i\not\mid\theta_{a} \ \ \ \mbox{for each term} \ a \ \mbox{of} \ TV^{-1} \ \mbox{and each} \ i\in [1,r]\setminus \mathcal{J},
\end{equation}
and that
\begin{equation}\label{equation fjnj(x)mid sigma(V)}
f_j^{n_j}\mid \theta_{\sigma(V)}  \ \ \ \mbox{ for each } j\in \mathcal{J}.
\end{equation}

For each term $a$ of $TV^{-1}$, let $\tilde{a}$ be the element of $\mathcal{S}_R$ such that
\begin{equation}\label{equation tilde a 1}
\theta_{\tilde{a}}\equiv \theta_{a}\pmod {f_i^{n_i}} \ \ \ \mbox{ for each } i\in [1,r]\setminus \mathcal{J}
\end{equation}
and
\begin{equation}\label{equation tilde a 2}
\theta_{\tilde{a}}\equiv 1_{\mathbb{F}_q}\pmod {f_j^{n_j}} \ \ \ \mbox{ for each } j \in\mathcal{J}.
\end{equation}
By \eqref{equation fi(x)notmid a}, \eqref{equation tilde a 1} and \eqref{equation tilde a 2}, we conclude that $\gcd(\theta_{\tilde{a}},f)=1_{\mathbb{F}_q}$, i.e.,
\begin{equation}\label{equation tilde a in U}
\tilde{a}\in U(\mathcal{S}_R)\ \ \ \mbox{for each term} \ a \ \mbox{of}\  TV^{-1}.
\end{equation}
By \eqref{equation fjnj(x)mid sigma(V)} and \eqref{equation tilde a 1}, we conclude that \begin{equation}\label{equation sigma(V)+tilde a}
\sigma(V)+\tilde{a}=\sigma(V)+a \ \ \ \mbox{for each term}\ \  a \ \mbox{of} \ TV^{-1}.
\end{equation}
By \eqref{equation length and D()} and \eqref{equation tilde a in U}, we have that $\prod\limits_{a\mid TV^{-1}}\tilde{a}$ is a nonempty sequence of elements in $U(\mathcal{S}_R)$ of length $|\prod\limits_{a\mid TV^{-1}}\tilde{a}|=|TV^{-1}|\geq {\rm D}(U(\mathcal{S}_R)\diagup K_t)$. It follows that there exists a {\bf nonempty} subsequence $$W\mid TV^{-1}$$  such that
$$\sigma(\prod\limits_{a\mid W}\tilde{a})\in K_t$$ which implies
\begin{equation}\label{equation sigma(tilde a)=0}
\sigma(V)+\sigma(\prod\limits_{a\mid W}\tilde{a})=\sigma(V).
\end{equation}
By \eqref{equation sigma(V)+tilde a} and \eqref{equation sigma(tilde a)=0}, we conclude that
$$\begin{array}{llll}
\sigma(T)&=& \sigma(TW^{-1}V^{-1})+(\sigma(V)+\sigma(W)) \\
&=& \sigma(TW^{-1}V^{-1})+(\sigma(V)+\sigma(\prod\limits_{a\mid W}\tilde{a}))\\
&=& \sigma(TW^{-1}V^{-1})+\sigma(V)\\
&=& \sigma(TW^{-1}),\\
\end{array}$$
and $T'=TW^{-1}$ is the desired proper subsequence of $T$.
This completes the proof of the theorem. \qed

\section{Concluding remarks}

We remark that if $R$ is the quotient ring of $\mathbb{F}_2[x]$, the conclusion ${\rm D}(\mathcal{S}_R)={\rm D}(U(\mathcal{S}_R))$ does not always hold true. For example, take $f=x*(x+1)*g\in \mathbb{F}_2[x]$ where $\gcd(x*(x+1),g)=1_{\mathbb{F}_2}$. Let $R=\mathbb{F}_2[x]\diagup (f)$. Take a sequence $T=a_1\cdot a_2\cdot \ldots\cdot a_{\ell}$, where $\theta_{a_1}=x$, $\theta_{a_2}=x+1$, and  $a_3\cdot\ldots\cdot a_{\ell}$ is a sequence of elements in $U(\mathcal{S}_R)$ of length $\ell-2={\rm D}(U(\mathcal{S}_R))-1$ which contains no nonempty subsequence $V$ with $\sigma(V)=0_{\mathcal{S}_R}$. It is easy to verify that $T$
is an irreducible sequence of length $\ell={\rm D}(U(\mathcal{S}_R))+1$, which implies that
${\rm D}(\mathcal{S}_R)\geq \ell+1={\rm D}(U(\mathcal{S}_R))+2$. Hence, we close this paper by proposing the following problem.

\medskip

\noindent {\bf Problem.} \ {\sl Let $R$ be a quotient ring of $\mathbb{F}_2[x]$ with $0\neq R\neq \mathbb{F}_2[x]$. Determine ${\rm D}(\mathcal{S}_R)-{\rm D}(U(\mathcal{S}_R))$.}

\bigskip

\noindent {\bf Acknowledgements}

\noindent  The author would like to thank the referee for his/her very useful suggestions.
This work is supported by NSFC (11301381, 11271207), Science and Technology Development Fund of Tianjin Higher
Institutions (20121003).

\end{document}